\documentclass{article}

\usepackage{microtype}
\usepackage{graphicx}
\usepackage{subfigure}
\usepackage{booktabs} 

\usepackage{hyperref}

\usepackage[accepted]{icml2019}

\icmltitlerunning{A Dynamical Systems Perspective on Nesterov Acceleration}
\usepackage[english,american]{babel}

\usepackage{url}
\usepackage{epsfig} 
\usepackage{times} 
\usepackage{amsmath}
\usepackage{amsfonts}
\usepackage{amssymb}  

\usepackage{xcolor}
\usepackage{transparent}
\usepackage{graphicx}
\graphicspath{{media/}}
\usepackage{tikz}

\usepackage{epstopdf}
\usepackage{units}
\usepackage{pgfplots}

\usepackage{pstool}

\graphicspath{{media/}}
\newcommand{\executeiffilenewer}[3]{%
 \ifnum\pdfstrcmp{\pdffilemoddate{#1}}%
 {\pdffilemoddate{#2}}>0%
 {\immediate\write18{#3}}\fi%
}

\usepackage{arydshln}
\usepackage{cleveref}
\usepackage{perpage}
\MakePerPage{footnote}

\newtheorem{theorem}{Theorem}[section]

\newtheorem{proposition}[theorem]{Proposition}

\newenvironment{proof}[1][Proof]{\begin{trivlist}
\item[\hskip \labelsep {\bfseries #1}]}{\end{trivlist}}

\newcommand{\qed}{\nobreak \ifvmode \relax \else
      \ifdim\lastskip<1.5em \hskip-\lastskip
      \hskip1.5em plus0em minus0.5em \fi \nobreak
      \vrule height0.75em width0.5em depth0.25em\fi}

\pgfplotsset{every tick label/.append style={font=\footnotesize}}

\begin{document}
\newcommand{\prox}{{\mathop{\mathrm{prox}}}}
\newcommand{\argmin}{\mathop{\mathrm{argmin}}}
\newcommand{\diag}{\mathop{\mathrm{diag}}}

\newcommand{\TT}{\ensuremath{\mathsf{\tiny{T}}}}
\newcommand{\T}{^{\TT}}
\newcommand{\Rot}[2][]{{#1}_{\textsc{\scriptsize{#2}}}}
\newcommand{\dRot}[2][]{\dot{#1}_{\textsc{\scriptsize{#2}}}}

\newcommand{\V}[2][]{{^{\textsc{\scriptsize{#2}}}{#1}}}
\newcommand{\dV}[2][]{{^{\textsc{\scriptsize{#2}}}\dot{#1}}}
\newcommand{\uV}[2][]{{_{\textsc{\scriptsize{#2}}}\vec{e}_{#1}}}

\newcommand{\tmat}[1]{\widetilde{#1}}
\newcommand{\tmatt}[2][]{{^{\textsc{\scriptsize{#2}}}\widetilde{#1}}}
\newcommand{\dtmatt}[2][]{{^{\textsc{\scriptsize{#2}}}\dot{\widetilde{#1}}}}

\newcommand{\sym}[1]{\mathrm{symm}(#1)}
\newcommand{\skw}[1]{\mathrm{skew}(#1)}

\newcommand{\diff}[1][]{\mathrm{d}#1}
\newcommand{\dt}{\diff t }

\newcommand{\Vct}[1]{\mathrm{vec}(#1)}

\newcommand{\lmax}{\lambda_{\mathrm{max}}}
\newcommand{\lmin}{\lambda_{\mathrm{min}}}

\newcommand{\tr}[1]{\mathrm{tr}(#1)}

\twocolumn[
\icmltitle{A Dynamical Systems Perspective on Nesterov Acceleration}



\icmlsetsymbol{equal}{*}

\begin{icmlauthorlist}
\icmlauthor{Michael Muehlebach}{berk}
\icmlauthor{Michael I. Jordan}{berk}
\end{icmlauthorlist}

\icmlaffiliation{berk}{Electrical Engineering and Computer Science Department, UC Berkeley, Berkeley, California, USA}

\icmlcorrespondingauthor{Michael Muehlebach}{michaelm@berkeley.edu}

\icmlkeywords{Optimization, Convex Optimization, Accelerated Gradient Method}

\vskip 0.3in
]



\printAffiliationsAndNotice{}  

\begin{abstract}
We present a dynamical system framework for understanding Nesterov's accelerated gradient method. In contrast to earlier work, our derivation does not rely on a vanishing step size argument. We show that Nesterov acceleration arises from discretizing an ordinary differential equation with a semi-implicit Euler integration scheme. We analyze both the underlying differential equation as well as the discretization to obtain insights into the phenomenon of acceleration. The analysis suggests that a curvature-dependent damping term lies at the heart of the phenomenon. We further establish connections between the discretized and the continuous-time dynamics.
\end{abstract}

\section{Introduction}
Many tasks in machine learning and statistics can be formulated as optimization problems. In its basic form, an optimization problem consists of finding a value $x^*$ that minimizes the function $f$, that is $f(x^*)\leq f(x)$ for all $x$, provided such a value exists. In the following, we consider the subclass of problems where the function $f$ is real-valued, convex, and has a Lipschitz-continuous gradient. These assumptions are often met, at least locally, in practical problem instances.

One of the most popular algorithms for computing $x^*$ up to a desired accuracy is the gradient method, due to its simplicity and due to the fact that it scales well to large problem sizes. It is also possible to obtain faster convergence rates within the family of gradient-based methods by considering accelerated gradient methods, which make use of two successive gradients \cite{NesterovOrgPaper}. 

Given the important role that acceleration has played not only in generating useful new algorithms but also in understanding natural limits on convergence rate, there have been many attempts to understand and characterize the phenomenon. \citet{Bubeck} suggest a modification of the accelerated gradient method that achieves the same convergence rate, but has a geometric interpretation. \citet{Allen-Zhu} present an algorithm that couples gradient and mirror descent, while achieving an accelerated convergence rate and \citet{LessardRecht} propose a general analysis framework, inspired by control theory, for comparing different first-order optimization algorithms including the accelerated gradient method. Other work includes \citet{Jelena}, who unify the analysis of first-order methods by imposing certain decay conditions, or \citet{Scieur}, where the accelerated gradient method is interpreted as a multi-step discretization of the gradient flow. 

An important step was taken by \citet{SuAcc} and \citet{KricheneAcc}, who showed that for a vanishing step size the trajectories of the accelerated gradient scheme approach the solutions of a certain second-order ordinary differential equation (ODE). The resulting ODE was analyzed in further detail by \citet{Attouch}.  It was placed within a variational framework by \citet{WibisonoVariational}, which led to the discovery of higher-order schemes that achieve even faster convergence rates.

In the current paper, we provide a different interpretation of the accelerated gradient method, one that is not based on beginning with the difference equation and deriving an underlying ODE via a vanishing step size argument. We go in the other direction, starting with a certain second-order ODE that has a clear interpretation as a mass-spring-damper system, and showing
that Nesterov acceleration can be derived by discretizing the ODE with a semi-implicit Euler method. In contrast to the continuous-time limit of the heavy-ball method \cite{PolyakHeavyBall}, our model includes a damping term that locally averages the curvature. This additional damping seems instrumental to the acceleration phenomenon. In addition to providing intuition, we believe that the model can be used to translate properties of the continuous-time system, which are, as we show, often easier to derive, to the resulting discrete-time algorithm. In particular, we show that under some regularity conditions fundamental geometric properties (phase-space area contraction and time-reversibility) are preserved by the discretization.

\emph{Notation and outline:} We focus our presentation on $n$-dimensional optimization problems over the real numbers, where the integer $n$ is finite. It will be clear that most derivations generalize in a straightforward way to non-Euclidean inner product spaces.  We also focus, again for simplicity of presentation, on strongly convex functions; that is, functions for which there exists a constant $\kappa\geq 1$ (the condition number) such that for any $\bar{x} \in \mathbb{R}^n$
\begin{equation*}
f(x)\geq f(\bar{x}) +\nabla f(\bar{x}) (x-\bar{x}) + \frac{L}{2 \kappa} |x-\bar{x}|^2, \forall x\in \mathbb{R}^n,
\end{equation*}
where $L$ is the Lipschitz constant of the gradient and $|\cdot|$ denotes the Euclidean norm. We also treat the non-strongly convex setting, however. The gradient and the Hessian (if it exists) evaluated at $x\in \mathbb{R}^n$ are denoted by $\nabla f(x) \in \mathbb{R}^n$ and $\Delta f(x) \in \mathbb{R}^{n\times n}$ respectively. We assume throughout the article that the function $f$ is smooth; i.e., that $f$ has a Lipschitz continuous gradient. We further assume that $f$ attains its minimum value at $x^*=0$ and that $f(0)=0$, which, in the case of a strongly convex function is without loss of generality.

The article is structured in the following way: Sec.~\ref{Sec:model} introduces the dynamical system model that serves as the foundation for our analysis. The dynamical system is shown to yield Nesterov's accelerated gradient method when discretized with the semi-implicit Euler scheme with a certain non-vanishing step size. In Sec.~\ref{Sec:Interpretation}, the continuous-time dynamics are analyzed and an intuitive interpretation as a mass-spring-damper system is provided. The dynamics are shown to have a curvature-dependent damping term, responsible for the acceleration. Sec.~\ref{Sec:Discretization} analyzes and motivates the discretization. It is shown that for a step size $T_\text{s}\in (0,1)$, the discrete-time dynamics preserve geometric properties of the continuous-time dynamics. Sec.~\ref{Sec:SimEx} presents a simulation example that illustrates the properties of the continuous-time dynamics and the discretization. The article concludes with a discussion in Sec.~\ref{Sec:conclusion}.
\section{Dynamical systems model}\label{Sec:model}
In this section we show that the accelerated gradient method presented in \citet{NesterovBook} (p.~81) results from a semi-implicit Euler discretization of the following ODE:
\begin{equation}
\ddot{x}(t)+2d \dot{x}(t) + \frac{1}{L \gamma^2} \nabla f(x(t)+\beta \dot{x}(t))=0,\label{eq:org}
\end{equation}
where $\gamma$ is a constant and
\begin{equation}
d:=\frac{1}{\sqrt{\kappa}+1}~\frac{1}{\gamma}, \quad \beta:=\frac{\sqrt{\kappa}-1}{\sqrt{\kappa}+1} ~\gamma.\label{eq:cons}
\end{equation}
Changing the constant $\gamma$ amounts to a rescaling of the solutions of \eqref{eq:org} in time. We set the constant $\gamma$ to one for ease of notation.

By the semi-implicit Euler scheme, we mean the following integration algorithm:
\begin{align}
q_{k+1}&=q_{k} + T_\text{s} \nabla T(p_{k+1}), \label{eq:Euler1}\\
p_{k+1}&=p_{k} + T_\text{s} (-\frac{1}{L}\nabla f(q_{k}) + f_{\text{NP}}(q_{k},p_{k})), \label{eq:Euler2}
\end{align}
where the corresponding continuous-time dynamics evolve according to
\begin{align}
\dot{q}(t)&=\nabla T(p(t)), \label{eq:ham1}\\
\dot{p}(t)&=-\frac{1}{L}\nabla f(q(t))+f_{\text{NP}}(q(t),p(t)), \label{eq:ham2}
\end{align}
and where the real-valued, continuously differentiable function $T(p)$ represents the kinetic energy, the continuous, real-valued function $f_\text{NP}(q,p)$ the non-potential forces, and the real number $T_\text{s}>0$ the step size. The total energy of the dynamic system described by \eqref{eq:ham1} and \eqref{eq:ham2} is given by
\begin{equation}
H(q,p):=T(p)+\frac{1}{L} f(q).
\end{equation}

The dynamic system \eqref{eq:org} can be brought to the form \eqref{eq:ham1} and \eqref{eq:ham2} by introducing
\begin{align}
T(p)&:=\frac{1}{2} |p|^2, \label{eq:T}\\
f_\text{NP}(q,p)&:=-2dp-\frac{1}{L}(\nabla f(q+\beta p)-\nabla f(q)), \label{eq:NP}
\end{align}
and identifying $q(t)$ with $x(t)$ and $p(t)$ with $\dot{x}(t)$. The discretization according to \eqref{eq:Euler1} and \eqref{eq:Euler2} with $T_\text{s}=1$ then yields
\begin{align}
p_{k+1}&=p_k - 2d p_k -\frac{1}{L} \nabla f(q_k+ \beta p_k)\nonumber\\
&=\beta p_k - \frac{1}{L}\nabla f(q_k+ \beta p_k), \label{eq:tmp1}\\
q_{k+1}&=q_k + p_{k+1} \nonumber\\
&=q_k + \beta p_k - \frac{1}{L} \nabla f(q_k + \beta p_k),
\end{align}
where the identity $2d + \beta = 1$ has been used in \eqref{eq:tmp1}. By defining $y_k:=q_k+\beta p_k$ and $x_k:=q_k$, we obtain
\begin{align}
x_{k+1}&=y_k - \frac{1}{L} \nabla f(y_k)\label{eq:nes1}\\
y_{k+1} &= x_{k+1} + \beta (x_{k+1} - x_{k})\label{eq:nes2},
\end{align}
which corresponds to the accelerated gradient scheme with constant step size \cite{NesterovBook} (p.~81).

A very similar derivation can be performed to obtain a formulation of the accelerated gradient method for smooth but non-strongly-convex functions. In particular, replacing $d$ and $\beta$ in \eqref{eq:org} with
\begin{equation}
\bar{d}(t):=\frac{3}{2(t+2)}, \quad \bar{\beta}(t):=\frac{t-1}{t+2}, \label{eq:dbar}
\end{equation}
and applying the above discretization with time-step $T_\text{s}=1$ yields the optimization algorithm that is used as a starting point in \citet{SuAcc}. The constants $d$, $\beta$, as well as the variables $\bar{d}(t)$, $\bar{\beta}(t)$ satisfy $2d+\beta=1$, respectively $2 \bar{d}(t)+ \bar{\beta}(t)=1$ for all $t\geq 0$.
\section{Interpretation of the model}\label{Sec:Interpretation}
In this section we discuss the interpretation of the ODE \eqref{eq:org} as a mass-spring-damper system with a curvature-dependent damping term. In addition, we show that the trajectories of \eqref{eq:org} converge exponentially with rate at least $1/(2\sqrt{\kappa})-1/(4\kappa)$, matching the well-known rate of the accelerated gradient method.

We start by reformulating the non-potential forces in \eqref{eq:NP} in the following way:\footnote{By assumption, the function $f$ has a Lipschitz-continuous gradient, hence the (Lebesgue) integral appearing in \eqref{eq:LebesgueIntegral} and \eqref{eq:defD} is well-defined even if the Hessian might not exist everywhere.}
\begin{align}
f_\text{NP}(q,p)&=-2dp - \frac{1}{L} \int_{0}^{\beta} \Delta f(q+\tau p) \text{d} \tau~p. \label{eq:LebesgueIntegral}
\end{align}
Thus, the dynamics \eqref{eq:org} can be rewritten as
\begin{equation}
\ddot{x}(t)+(2dI+D_{x\dot{x}})~\dot{x}(t) +\frac{1}{L} \nabla f(x(t))=0,
\end{equation}
where $I\in \mathbb{R}^{n\times n}$ denotes identity matrix and where the damping $D_{x\dot{x}}$ is defined as
\begin{equation}
D_{x\dot{x}}:=\frac{1}{L} \int_{0}^{\beta} \Delta f(x(t) + \tau \dot{x}(t))~\text{d}\tau.\label{eq:defD}
\end{equation}
To lighten notation the dependence of $D_{x\dot{x}}$ on $x(t)$ and $\dot{x}(t)$ is indicated with subscripts.
Direct calculations show that
\begin{align}
\frac{\text{d}}{\text{d}t} H(x(t),\dot{x}(t))&=f_\text{NP}(x(t),\dot{x}(t))\T \dot{x}(t) \nonumber\\
&=-\dot{x}(t)\T (2dI+D_{x\dot{x}}) \dot{x}(t),
\end{align}
which bounds the energy dissipation by
\begin{equation}
-|\dot{x}(t)|^2 \leq \frac{\text{d}}{\text{d}t} H(x(t),\dot{x}(t)) \leq -(2d +\beta/\kappa) |\dot{x}(t)|^2.
\end{equation}
Note that the identity $2d + \beta =1$, as well as the assumption of $f$ being smooth and strongly convex have been used for obtaining these inequalities. In the non-strongly convex case the bounds reduce to 
\begin{equation}
-|\dot{x}(t)|^2 \leq \frac{\text{d}}{\text{d}t} H(x(t),\dot{x}(t)) \leq -2\bar{d}(t) |\dot{x}(t)|^2,
\end{equation}
for $t\geq 1$, respectively
\begin{equation}
-2\bar{d}(t) |\dot{x}(t)|^2 \leq \frac{\text{d}}{\text{d}t} H(x(t),\dot{x}(t)) \leq -|\dot{x}(t)|^2,
\end{equation}
for $t\in (0,1)$.

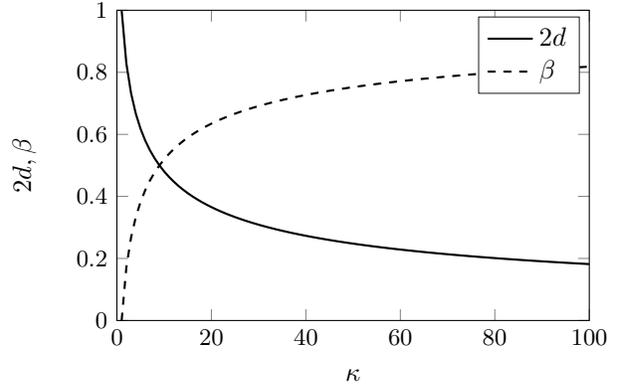
\begin{figure}
\newlength{\figurewidth}
\newlength{\figureheight}
\setlength{\figurewidth}{.8\columnwidth}
\setlength{\figureheight}{.5\columnwidth}
%
%
\begin{tikzpicture}

\begin{axis}[%
width=0.95092\figurewidth,
height=\figureheight,
at={(0\figurewidth,0\figureheight)},
scale only axis,
xmin=0,
xmax=100,
xlabel={$\kappa$},
ymin=0,
ymax=1,
ylabel={$2d, \beta$},
legend style={legend cell align=left,align=left,draw=white!15!black}
]
\addplot [color=black,solid,thick]
  table[row sep=crcr]{%
1	1\\
2	0.82842712474619\\
3	0.732050807568877\\
4	0.666666666666667\\
5	0.618033988749895\\
6	0.579795897113271\\
7	0.548583770354863\\
8	0.522407749927483\\
9	0.5\\
10	0.480506146704084\\
11	0.46332495807108\\
12	0.448018475479592\\
13	0.434258545910665\\
14	0.421793444119068\\
15	0.410426192315345\\
16	0.4\\
17	0.390388203202208\\
18	0.381487139661092\\
19	0.373210993726741\\
20	0.365487995263114\\
21	0.358257569495584\\
22	0.351468167602231\\
23	0.345075593028429\\
24	0.339041694397074\\
25	0.333333333333333\\
26	0.327921561087423\\
27	0.322780955592818\\
28	0.31788908312068\\
29	0.313226057652465\\
30	0.30877417758977\\
31	0.304517624188668\\
32	0.30044220964467\\
33	0.296535165408627\\
34	0.292784963323958\\
35	0.289181163711742\\
36	0.285714285714286\\
37	0.282375696127679\\
38	0.279157513673999\\
39	0.276052526231495\\
40	0.273054118991629\\
41	0.270156211871642\\
42	0.267353204800383\\
43	0.264639929728667\\
44	0.262011608405153\\
45	0.259463815113608\\
46	0.256992443694456\\
47	0.254593678278306\\
48	0.252263967245766\\
49	0.25\\
50	0.247798686198591\\
51	0.245657137141714\\
52	0.243572649055999\\
53	0.241542688049251\\
54	0.239564876541492\\
55	0.237636981003543\\
56	0.235756900856287\\
57	0.233922658402527\\
58	0.232132389679435\\
59	0.2303843361334\\
60	0.228676837031011\\
61	0.227008322530222\\
62	0.22537730734465\\
63	0.223782384941735\\
64	0.222222222222222\\
65	0.22069555463433\\
66	0.219201181681106\\
67	0.217737962784014\\
68	0.216304813469711\\
69	0.214900701850532\\
70	0.213524645372196\\
71	0.212175707805039\\
72	0.210852996457425\\
73	0.209555659592154\\
74	0.208282884028565\\
75	0.207033892914713\\
76	0.205807943655503\\
77	0.204604325984003\\
78	0.20342236016436\\
79	0.202261395315784\\
80	0.20112080784808\\
81	0.2\\
82	0.198898398472529\\
83	0.197815453149861\\
84	0.196750635901486\\
85	0.195703439459354\\
86	0.194673376364605\\
87	0.19365997797881\\
88	0.1926627935551\\
89	0.191681389364923\\
90	0.19071534787652\\
91	0.189764266981543\\
92	0.188827759266493\\
93	0.187905451325934\\
94	0.186996983114681\\
95	0.186102007336361\\
96	0.185220188865952\\
97	0.184351204204086\\
98	0.183494740961065\\
99	0.182650497368698\\
100	0.181818181818182\\
};
\addlegendentry{$2d$};

\addplot [color=black,dashed,thick]
  table[row sep=crcr]{%
1	0\\
2	0.17157287525381\\
3	0.267949192431123\\
4	0.333333333333333\\
5	0.381966011250105\\
6	0.420204102886729\\
7	0.451416229645137\\
8	0.477592250072517\\
9	0.5\\
10	0.519493853295916\\
11	0.53667504192892\\
12	0.551981524520408\\
13	0.565741454089335\\
14	0.578206555880932\\
15	0.589573807684655\\
16	0.6\\
17	0.609611796797792\\
18	0.618512860338908\\
19	0.626789006273259\\
20	0.634512004736886\\
21	0.641742430504416\\
22	0.648531832397769\\
23	0.654924406971571\\
24	0.660958305602926\\
25	0.666666666666667\\
26	0.672078438912577\\
27	0.677219044407182\\
28	0.68211091687932\\
29	0.686773942347535\\
30	0.69122582241023\\
31	0.695482375811332\\
32	0.69955779035533\\
33	0.703464834591373\\
34	0.707215036676042\\
35	0.710818836288258\\
36	0.714285714285714\\
37	0.717624303872321\\
38	0.720842486326001\\
39	0.723947473768505\\
40	0.726945881008371\\
41	0.729843788128358\\
42	0.732646795199617\\
43	0.735360070271333\\
44	0.737988391594847\\
45	0.740536184886392\\
46	0.743007556305544\\
47	0.745406321721694\\
48	0.747736032754234\\
49	0.75\\
50	0.752201313801409\\
51	0.754342862858286\\
52	0.756427350944001\\
53	0.758457311950749\\
54	0.760435123458508\\
55	0.762363018996457\\
56	0.764243099143713\\
57	0.766077341597473\\
58	0.767867610320565\\
59	0.7696156638666\\
60	0.771323162968989\\
61	0.772991677469778\\
62	0.77462269265535\\
63	0.776217615058265\\
64	0.777777777777778\\
65	0.77930444536567\\
66	0.780798818318894\\
67	0.782262037215986\\
68	0.783695186530289\\
69	0.785099298149468\\
70	0.786475354627804\\
71	0.787824292194961\\
72	0.789147003542575\\
73	0.790444340407846\\
74	0.791717115971435\\
75	0.792966107085287\\
76	0.794192056344497\\
77	0.795395674015997\\
78	0.79657763983564\\
79	0.797738604684216\\
80	0.79887919215192\\
81	0.8\\
82	0.801101601527471\\
83	0.802184546850139\\
84	0.803249364098514\\
85	0.804296560540646\\
86	0.805326623635395\\
87	0.80634002202119\\
88	0.8073372064449\\
89	0.808318610635077\\
90	0.80928465212348\\
91	0.810235733018457\\
92	0.811172240733507\\
93	0.812094548674066\\
94	0.813003016885319\\
95	0.813897992663639\\
96	0.814779811134048\\
97	0.815648795795914\\
98	0.816505259038935\\
99	0.817349502631302\\
100	0.818181818181818\\
};
\addlegendentry{$\beta$};

\end{axis}
\end{tikzpicture}%
\caption{The figure shows the constants $2d$ and $\beta$ as a function of the maximum curvature $\kappa$. The constant $2d$ decreases with $1/\sqrt{\kappa}$, whereas $\beta$ increases with $1/\sqrt{\kappa}$ (for large values of $\kappa$).}
\label{Fig:2dbeta}
\end{figure}

Summarizing, it may therefore be concluded (the conclusions in the non-strongly-convex case are analogous):
\begin{itemize}
\item In the absence of damping, i.e., $d=D_{x\dot{x}}\equiv 0$, the total energy $H(x(t),\dot{x}(t))$ is conserved. The trajectories are thus confined to the level sets of $H(x(t),\dot{x}(t))$.
\item In the presence of damping, energy is dissipated. As a result, La Salle's theorem~\citep[see, for example,][Ch.~5.4]{SastryNonlinear} implies that the trajectories of \eqref{eq:org} converge asymptotically to the origin. Similar conclusions can be drawn in the non-convex case.
\item The damping coefficient $D_{x\dot{x}}$ can be interpreted as a local average of the curvature near $x(t)$. The amount of averaging depends on the maximum curvature of the objective function and on the current velocity. More averaging is performed at higher velocities (larger $\dot{x}(t)$) and for a larger maximum curvature (large value of $\kappa$). The latter is due to the fact that the coefficient $\beta$ increases monotonically with $\kappa$.
\item The total damping is a linear combination of constant damping and local curvature-dependent damping. The two sources of damping are balanced by the coefficient $2d$ (constant damping) that decreases with larger $\kappa$ and the coefficient $\beta$ (local curvature-dependent damping) that increases with larger $\kappa$. The total amount of damping remains constant in the sense that $2d+\beta=1$ for any maximum curvature $\kappa$. The constants $2d$ and $\beta$ are shown in Fig.~\ref{Fig:2dbeta} for various values of $\kappa$.
\item Energy is dissipated even for $\kappa\rightarrow \infty$ due to the local curvature-dependent damping.
\end{itemize}

The following two propositions provide an estimate of the convergence rate of the differential equation \eqref{eq:org}, the proof of which can be found in the appendix, App.~\ref{App:ProofSC} and App.~\ref{App:ProofNSC}.

\begin{proposition}
(Strongly convex case) The trajectories of \eqref{eq:org} decay exponentially to the origin with rate at least $1/(2\sqrt{\kappa})-1/(4\kappa)$. \label{Prop:DecayExp}
\end{proposition}

\begin{proposition}
(Non-strongly-convex case) When replacing $d$ with $\bar{d}(t)$ and $\beta$ with $\bar{\beta}(t)$ in \eqref{eq:org}, $f(x(t))$ decays at least with $\mathcal{O}(1/t^2)$. \label{Prop:DecaySmooth}
\end{proposition}

We conclude the section with the following remarks about Prop.~\ref{Prop:DecayExp} and Prop.~\ref{Prop:DecaySmooth}:
\begin{itemize}
\item Prop.~\ref{Prop:DecayExp} and Prop.~\ref{Prop:DecaySmooth} are based on Lyapunov functions of the type
\begin{equation}
\frac{1}{2} |a q + p |^2 + \frac{1}{L} f(q), \label{eq:Lyap}
\end{equation}
where $a$ is related to the rate of convergence and is either constant, in the strongly convex case, or time-varying, in the non-strongly-convex case. 
\item In both cases, the choice of $a$ is motivated by convenience, such that the resulting cross-term $q\T p$ vanishes for the given parameters $d$ and $\beta$ (or $\bar{d}$ and $\bar{\beta}$ in the non-strongly-convex case). However, the proof can be extended to derive more general sufficiency conditions on the parameters $d$ and $\beta$ (or $\bar{d}$ and $\bar{\beta}$) that provide an accelerated rate.
\item The asymptotic convergence rate of $\mathcal{O}(1/\sqrt{\kappa})$ in Prop.~\ref{Prop:DecayExp} can also be established using a Lyapunov function that is inspired by the Popov criterion; see for example \cite{SastryNonlinear} (Ch.~6.2.2). The Popov criterion gives a sufficient condition for the closed-loop stability of a linear time-invariant system subject to a static nonlinear feedback path, and is applicable to the ODE \eqref{eq:org}.
\end{itemize}
\section{The discretization}\label{Sec:Discretization}
The following section highlights that fundamental geometric properties, such as the phase-space area contraction rate and the time-reversibility of the dynamics \eqref{eq:org} are preserved through the discretization. In addition to providing insights, both properties will be used to bound the convergence rate in a certain worst-case sense. The section concludes with a proposition stating that the discretized dynamics indeed converge at the given accelerated rate. Unfortunately, the proof relies on a Lyapunov function, similar to \eqref{eq:Lyap}, and does not seem to provide additional insights into the discretization process. For simplifying the presentation, the result concerning the phase-space area contraction is stated for a one-dimensional objective function, while the more general case can be found in App.~\ref{App:ProofContractionIII}. The map $(q_k,p_k)\rightarrow (q_{k+1},p_{k+1})$ given by the discretized dynamics \eqref{eq:Euler1} and \eqref{eq:Euler2} is denoted by $\psi$.

The discretization \eqref{eq:Euler1} and \eqref{eq:Euler2} can be divided into two parts, a non-conservative step that involves an update of the momentum based on the non-potential forces $f_{\text{NP}}$, 
\begin{align}
\bar{q}_{k+1}&=q_{k}, \quad \bar{p}_{k+1}=p_{k} + T_{\text{s}} (f_{\text{NP}}(q_{k},p_{k})), \label{eq:contract}
\end{align}
followed by a symplectic Euler step based on the conservative part of the system,
\begin{align}
q_{k+1} &= \bar{q}_{k+1}+T_{\text{s}} \nabla T(p_{k+1}), \label{eq:SympEuler1}\\
p_{k+1} &= \bar{p}_{k+1}+T_{\text{s}} (-\frac{1}{L} \nabla f(\bar{q}_{k+1})). \label{eq:SympEuler2}
\end{align}
The symplectic Euler scheme is a well-known integration scheme for Hamiltonian systems~\citep[see, for example,][Ch.~VI.3]{HairerGeom}. It is one of the simplest symplectic integration schemes and is known to be energy consistent (nearly energy conserving) over exponentially long time intervals~\cite{HairerGeom}.

The combination of \eqref{eq:contract} with \eqref{eq:SympEuler1} and \eqref{eq:SympEuler2} leads to a phase-space area contraction, which we quantify with the following proposition. (Here formulated for $n=1$.)

\begin{proposition}
Let $\partial \Gamma_k$ be a simple closed contour in $\mathbb{R}^2$ and let the signed area corresponding to $\partial \Gamma_k$ be defined as
\begin{equation}
A_k:=\int_{\Gamma_k} \text{d}q_k\wedge \text{d}p_k,
\end{equation}
where $\Gamma_k$ describes all points enclosed by $\partial \Gamma_k$.\footnote{The sign is inferred from the contour, that is, a counter-clockwise direction yields a positive sign: $\text{d}q_k \wedge \text{d}p_k=\text{d}q_k \text{d}p_k$, a clockwise direction yields a negative sign, i.e. $\text{d}q_k \wedge \text{d}p_k=-\text{d}q_k \text{d}p_k$, where $\text{d}q_k \text{d}p_k$ refers to the standard area measure.}
The integration scheme \eqref{eq:Euler1} and \eqref{eq:Euler2} maps $\partial \Gamma_k$ to $\partial \Gamma_{k+1}$ such that 
\begin{align}
A_{k+1}-A_{k}&=-T_\text{s} \int_{\partial \Gamma_{k}} f_\text{NP}(q_k,p_k) \text{d}q_k\label{eq:proofeq11d}\\
&=-T_\text{s} \int_{\Gamma_k} 2d + \frac{\beta}{L} \Delta f(q_k+\beta p_k)~\text{d}q_k\wedge\text{d}p_k.\label{eq:proofeq21d}
\end{align}
\label{Prop:Contraction}
\end{proposition}

\begin{figure}
    \centering
    \def\svgwidth{1.3\columnwidth}
    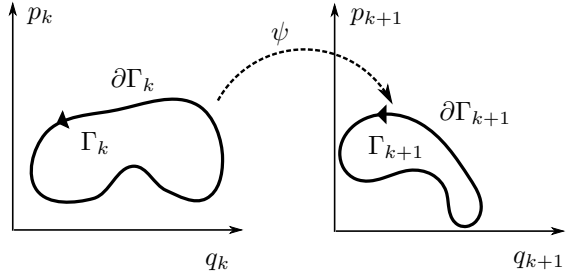
    \vspace*{-10mm}
    \caption{The figure illustrates the setting of Prop.~\ref{Prop:Contraction}. The direction of the contours is indicated by the two arrows, and as a result, both $A_k$, enclosed by $\partial \Gamma_k$, and $A_{k+1}$, enclosed by $\partial \Gamma_{k+1}$, have a positive sign. The proposition states that the areas $A_k$, $k=0,1,\dots$, contract with a rate matching the continuous-time dynamics for $T_\text{s}\in (0,1)$.}
    \label{Fig:Contraction}
\end{figure}

The setting of the proposition is illustrated by Fig.~\ref{Fig:Contraction} and its proof can be found in App.~\ref{App:Gen}.

Provided that the area $A_k$ has a positive sign (for a negative sign the inequalities are reversed), the integral on the right-hand side of \eqref{eq:proofeq21d} is bounded by
\begin{multline}
(2d + \frac{\beta}{\kappa}) A_k \leq \int_{\Gamma_{k}} \hspace{-2pt}2d + \frac{\beta}{L} \Delta f(q_k+\beta p_k) \text{d}q_k\text{d}p_k \leq A_k,\label{eq:boundA}
\end{multline}
which follows from the fact that $f$ is strongly convex (yielding the lower bound) and smooth (yielding the upper bound), and the identity $2d+\beta=1$. This concludes that the phase-space area contracts for $T_\text{s}\in (0,2)$. In case $T_\text{s} \in (0,1]$ the area does not change sign; i.e., the contour is guaranteed to preserve its original orientation. For $T_\text{s}=1$, which yields the accelerated gradient method, the contour might contract to a single point.

With very similar arguments it can be shown that the phase-space area contraction rate of the continuous-time dynamics \eqref{eq:org} is given by
\begin{align}
\dot{A}(t)&=-\int_{\partial \Gamma(t)} f_\text{NP}(q,p) \text{d}q\\
&=-\int_{\Gamma(t)} 2d + \frac{\beta}{L} \Delta f(q+\beta p)~\text{d}q\wedge\text{d}p, \label{eq:ctac}
\end{align}
where $\partial \Gamma(t)$ denotes a simple closed contour in $\mathbb{R}^2$ that evolves according to the dynamics \eqref{eq:org}, and $A(t)$ is the signed area enclosed by the contour. It can therefore be concluded that the phase-space area of the discretized dynamics essentially contracts at the same rate as the continuous-time counterpart. More precisely, we have:
\begin{proposition}
The contraction of the signed phase-space area enclosed by an energy level-set with the discretized dynamics is larger than the contraction with the continuous-time dynamics over the length of one time step. \label{Prop:ContractionIII}
\end{proposition}

The proof relies on the mean-value theorem and can be found in App.~\ref{App:ProofContractionIII}. 

We show next that the discrete integration describes a homeomorphism (a continuous bijection that has a continuous inverse) for $T_\text{s}\in (0,1)$. This implies that the discrete-time dynamics are time reversible, which parallels the continuous-time counterpart.

\begin{proposition}\label{Prop:Hom}
The discrete integration given by \eqref{eq:Euler1} and \eqref{eq:Euler2} describes a homeomorphism provided that $T_\text{s}\in (0,1)$.
\end{proposition}

The proof can be found in App.~\ref{App:ProofOfHom}. 

Combining the above proposition with the phase-space area contraction yields the following result:
\begin{proposition}
For every $T_\text{s} \in (0,1)$ there exists at least one trajectory that converges linearly with rate $1-T_\text{s}\mathcal{O}(1/\sqrt{\kappa})$.
\label{Prop:ContractionII}
\end{proposition}
\begin{proof}
We consider the closed contour $\partial \Gamma_0$, given by a circle with radius $R>0$ centered at the origin (with positive orientation). The map $\psi$ is a homeomorphism, and therefore, as the contour evolves, the origin will always be contained in its interior. The area $A_0$ enclosed by the contour at time zero is $R^2\pi$. We claim that at time $k>0$, there is at least one trajectory that is at a distance at least
\begin{equation}
R_k:=R (1-T_\text{s}(2d+\beta/\kappa))^{k/2}
\end{equation}
from the origin. For the sake of contradiction we assume that the claim is false. This implies, however, that all trajectories starting from $\partial \Gamma_0$ remain outside a closed ball of radius $R_k$ centered at the origin, and hence $A_k>R_k^2 \pi$. This is contradicting the fact that $A_k$ decays at least with $(1-T_\text{s} (2d+\beta/\kappa))^k$, according to Prop.~\ref{Prop:Contraction} and \eqref{eq:boundA}, which concludes the proof. \qed
\end{proof}

We show in the appendix (App.~\ref{App:HardConvergence}) that all trajectories indeed converge at an accelerated rate, as suggested by Prop.~\ref{Prop:ContractionII}:
\begin{proposition}
For every $T_\text{s}\in (0,1]$ the trajectories of the discrete-time dynamics \eqref{eq:Euler1} and \eqref{eq:Euler2} converge linearly with rate at least $1-T_s\mathcal{O}(1/\sqrt{\kappa})$. \label{Prop:HardConvergence}
\end{proposition}

\begin{figure*}
\setlength{\figurewidth}{2\columnwidth}
\setlength{\figureheight}{0.8\columnwidth}
\centering
\setlength\fboxsep{0pt}
\setlength\fboxrule{0pt}
 \tikzset{mark size=0.5pt}
\fbox{\input{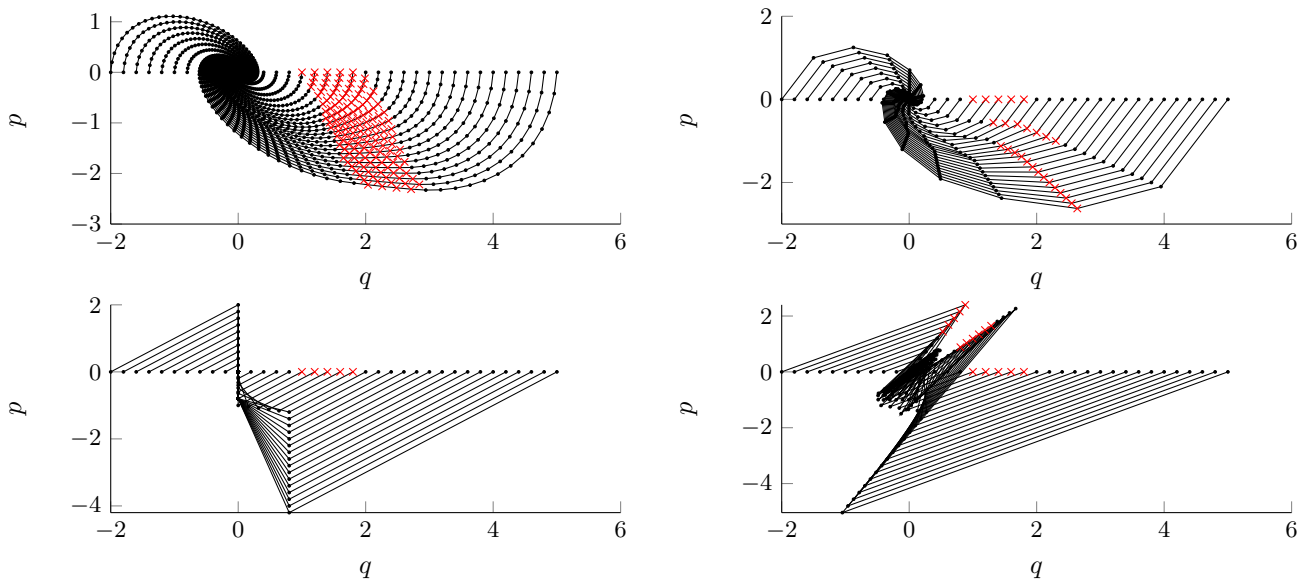}}
\vspace*{-5mm}
\caption{The plot shows the evolution of different trajectories in the phase space. The trajectories are obtained by applying the discrete integration algorithm \eqref{eq:Euler1} and \eqref{eq:Euler2} with the time steps $T_\text{s}=0.1,0.5,1,1.2$ (top left, top right, bottom left, bottom right). The points $(q,p)$ with $q+\beta p\in [1,2)$ are marked with red crosses.}
\label{Fig:PhaseSpace}
\end{figure*}

\section{Simulation example}\label{Sec:SimEx}
This section presents a numerical example that illustrates certain properties of the discrete integration. 

\begin{figure}
\setlength{\figurewidth}{.8\columnwidth}
\setlength{\figureheight}{.45\columnwidth}
%
%
\begin{tikzpicture}

\begin{axis}[%
width=0.95092\figurewidth,
height=\figureheight,
at={(0\figurewidth,0\figureheight)},
scale only axis,
xmin=-1,
xmax=4,
xlabel={$x$},
ymin=-5,
ymax=30,
ylabel={$f(x)$, $\nabla f(x)$},
legend style={legend cell align=left,align=left,draw=white!15!black},
legend pos={north west}
]
\addplot [color=black,solid]
  table[row sep=crcr]{%
-1	2.5\\
-0.9	2.025\\
-0.8	1.6\\
-0.7	1.225\\
-0.6	0.9\\
-0.5	0.625\\
-0.4	0.4\\
-0.3	0.225\\
-0.2	0.1\\
-0.1	0.025\\
0	0\\
0.1	0.025\\
0.2	0.1\\
0.3	0.225\\
0.4	0.4\\
0.5	0.625\\
0.6	0.9\\
0.7	1.225\\
0.8	1.6\\
0.9	2.025\\
1	2.5\\
1.1	3.005\\
1.2	3.52\\
1.3	4.045\\
1.4	4.58\\
1.5	5.125\\
1.6	5.68\\
1.7	6.245\\
1.8	6.82\\
1.9	7.405\\
2	8\\
2.1	8.625\\
2.2	9.3\\
2.3	10.025\\
2.4	10.8\\
2.5	11.625\\
2.6	12.5\\
2.7	13.425\\
2.8	14.4\\
2.9	15.425\\
3	16.5\\
3.1	17.625\\
3.2	18.8\\
3.3	20.025\\
3.4	21.3\\
3.5	22.625\\
3.6	24\\
3.7	25.425\\
3.8	26.9\\
3.9	28.425\\
4	30\\
};
\addlegendentry{$f(x)$};

\addplot [color=red,solid]
  table[row sep=crcr]{%
-1	-5\\
-0.9	-4.5\\
-0.8	-4\\
-0.7	-3.5\\
-0.6	-3\\
-0.5	-2.5\\
-0.4	-2\\
-0.3	-1.5\\
-0.2	-1\\
-0.1	-0.5\\
0	0\\
0.1	0.5\\
0.2	1\\
0.3	1.5\\
0.4	2\\
0.5	2.5\\
0.6	3\\
0.7	3.5\\
0.8	4\\
0.9	4.5\\
1	5\\
1.1	5.1\\
1.2	5.2\\
1.3	5.3\\
1.4	5.4\\
1.5	5.5\\
1.6	5.6\\
1.7	5.7\\
1.8	5.8\\
1.9	5.9\\
2	6\\
2.1	6.5\\
2.2	7\\
2.3	7.5\\
2.4	8\\
2.5	8.5\\
2.6	9\\
2.7	9.5\\
2.8	10\\
2.9	10.5\\
3	11\\
3.1	11.5\\
3.2	12\\
3.3	12.5\\
3.4	13\\
3.5	13.5\\
3.6	14\\
3.7	14.5\\
3.8	15\\
3.9	15.5\\
4	16\\
};
\addlegendentry{$\nabla f(x)$};

\end{axis}
\end{tikzpicture}%
\label{Fig:plotF}
\caption{The plot shows the function $f$ and its gradient $\nabla f$ used for the numerical example.}
\end{figure}
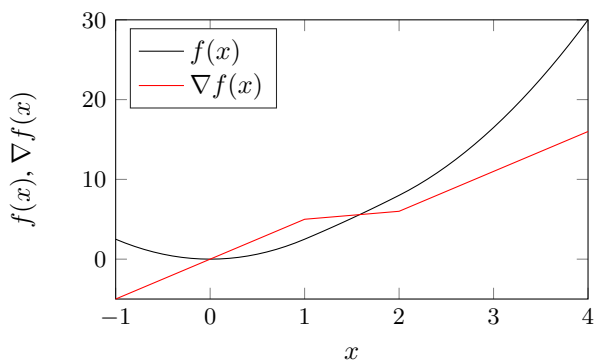

We choose a function $f$ with the following gradient
\begin{equation}
\nabla f(x)=\begin{cases} \kappa x, \quad &x<1,\\ \kappa-1+x, \quad &1\leq x <2,\\ 1-\kappa+\kappa x, \quad &x\geq 2, \end{cases} 
\end{equation}
where the condition number $\kappa$ is set to 5.
The function and its gradient are shown in Fig.~\ref{Fig:plotF}, and is inspired by the counterexample provided in \citet{LessardRecht}. The integration algorithm given by \eqref{eq:Euler1} and \eqref{eq:Euler2} is applied to the initial conditions $(q_0,0)$, where $q_0$ is varied from $-2$ to $5$ in steps of $0.2$. The step size $T_\text{s}$ is successively increased from $0.1$ to $1.2$. The evolution of the corresponding trajectories in the phase space is shown in Fig.~\ref{Fig:PhaseSpace}. For a time step of $T_\text{s}=0.1$ the trajectories resemble the solutions of \eqref{eq:org}. As the time step is increased the trajectories are more and more ``folded" in clockwise direction. For the time step $T_\text{s}=1$, which corresponds to the accelerated gradient method, multiple initial conditions are mapped to the same point (for example the point $q=0$, $p=-0.8$) due to the course of the integration. The trajectories starting from an initial condition $q_0<1$, $p_0=0$ converge exactly in two steps. If the time step is increased above $T_\text{s}=1$ the discrete integration map $\psi$ is no longer orientation preserving, and the resulting motion is much less structured, as shown in Fig.~\ref{Fig:PhaseSpace} (bottom right). For a time step of $T_\text{s}=1.3$ trajectories starting from $(q_0,0)$ with $q_0\geq 4.4$ are found to diverge.

\section{Conclusions}\label{Sec:conclusion}
We have shown that Nesterov's accelerated gradient method can be interpreted as the discretization of a certain ODE with the semi-implicit Euler scheme. The differential equation describes a mass-spring-damper system with a curvature-dependent damping term, which seems essential for the acceleration. The analysis suggests that geometric properties (time reversibility and phase-space area contraction) are conserved by the discretization. 

The discretization can be divided into two steps, a contraction step based on the non-conservative part and a symplectic Euler step based on the conservative part of the dynamics. The fact that the symplectic Euler scheme is long-term energy consistent makes the discretization appear very natural. However, as the steps do not commute, a straightforward application of backwards error analysis, which would rigorously justify this intuition, seems difficult. 

Even though we characterized and related the phase-space area contraction property of the discretized dynamics to their continuous-time counterpart, the area contraction property alone does not seem to be enough for deriving a convergence rate of $f(q)$.

\section*{Acknowledgements}
We thank the Branco Weiss Fellowship, administered by ETH Zurich, for the generous support and the
Office of Naval Research under grant number N00014-18-1-2764.

\bibliography{literature}
\bibliographystyle{icml2019}

\newpage
\appendix
\section{Appendix}
\subsection{Proof of Prop.~\ref{Prop:DecayExp}} \label{App:ProofSC}
We claim that the following function is a Lyapunov function for the differential equation \eqref{eq:org} (with $q(t)=x(t)$, $p(t)=\dot{x}(t)$):
\begin{equation}
V(q,p):=\frac{1}{2}|a q+p|^2 + \frac{1}{L} f(q),
\end{equation}
where 
\begin{equation}
a:=d+\frac{\beta}{2\kappa}=\frac{1}{\sqrt{\kappa}}-\frac{1}{2\kappa}.
\end{equation} 
The value of the constant $a$ is motivated below. It can be verified that $V$ is positive definite. We will omit the dependence on time to simplify notation. Evaluating the time derivative of $V$ along the trajectories of \eqref{eq:org} yields
\begin{multline}
\dot{V}=-\frac{a}{L} q\T \nabla f(y) + a(a-2d) q\T p + (a-2d)|p|^2 \\
- \frac{1}{L} p\T (\nabla f(y)-\nabla f(q)), \label{eq:tmp10}
\end{multline}
with $y:=q + \beta p$. Due to the fact that $f$ is smooth and strongly convex, the following two inequalities hold:
\begin{align}
-y\T \nabla f(y) \leq - f(y) - \frac{L}{2\kappa} |y|^2, \label{eq:tmp11}\\
f(y) \geq f(q) + \nabla f(q)\T \beta p + \frac{L}{2\kappa}\beta^2 |p|^2. \label{eq:tmp12}
\end{align}
These can be combined to
\begin{equation*}
-y\T \nabla f(y) \leq -f(q) - \nabla f(q)\T \beta p -\frac{L}{2 \kappa}(|y|^2 + \beta^2 |p|^2),
\end{equation*}
which implies, by definition of $y$, that
\begin{multline}
-q\T \nabla f(y) \leq -f(q) - \frac{L}{2 \kappa} (|q|^2 + 2 \beta q\T p + 2 \beta^2 |p|^2)\\
+\beta p\T (\nabla f(y)-\nabla f(q)). \label{eq:tmp13}
\end{multline}
Substituting the above bound on $-q\T \nabla f(y)$ in \eqref{eq:tmp10} and rearranging terms yields
\begin{multline}
\dot{V}=-aV + \frac{a}{2} \left(a^2-\frac{1}{\kappa}\right) |q|^2 + \left(\frac{3a}{2}-2d -\frac{a \beta^2}{\kappa}\right) |p|^2\\
+a\left( 2a-2d -\frac{\beta}{\kappa}\right) q\T p\\
-\frac{1}{L} (1-\beta a) p\T (\nabla f(y)-\nabla f(q)). \label{eq:tmp14}
\end{multline}
The constant $a$ is deliberately chosen such that the cross-term $q\T p$ in \eqref{eq:tmp14} vanishes. Due to the fact that $1-\beta a\geq 0$ and $\beta\geq 0$, the last term of \eqref{eq:tmp14} can be bounded by
\begin{equation*}
-\frac{1}{L} (1-\beta a) p\T (\nabla f(y) - \nabla f(q)) \leq -(1-\beta a)\frac{\beta}{\kappa} |p|^2,
\end{equation*}
which implies,
\begin{multline}
\dot{V}\leq-aV + \frac{a}{2} \left(a^2-\frac{1}{\kappa}\right) |q|^2 + \left(\frac{3a}{2} -2d -\frac{\beta}{\kappa}\right) |p|^2.
\end{multline}
It can be verified that $a^2-1/\kappa\leq 0$ and $(3a/2-2d-\beta/\kappa)\leq 0$ for all $\kappa \geq 1$. Hence, $\dot{V}\leq -a V$, which, by applying Gr\"{o}nwall's inequality, establishes the proposition. \qed

\subsection{Proof of Prop.~\ref{Prop:DecaySmooth}}\label{App:ProofNSC}
The proof is analogous to the proof of Prop.~\ref{Prop:DecayExp}, and hinges on the following Lyapunov-like function
\begin{equation}
\bar{V}(t):=\frac{1}{2}|\bar{a}(t) q(t)+p(t)|^2 + \frac{1}{L} f(q(t)),
\end{equation}
where $\bar{a}(t):=2/(t+2)$, and $q(t)=x(t), p(t)=\dot{x}(t)$ are the solutions of \eqref{eq:org}, when $d$ is replaced with $\bar{d}(t)$ and $\beta$ with $\bar{\beta}(t)$. The dependence on time will be omitted to simplify notation. Evaluating the time derivative of $\bar{V}$, applying inequalities \eqref{eq:tmp11} and \eqref{eq:tmp12} (with $\kappa\rightarrow \infty$) as in App.~\ref{App:ProofSC}, and rearranging terms yields
\begin{multline}
\dot{\bar{V}}\leq -\bar{a} \bar{V} + \left(\frac{\bar{a}^3}{2} + \bar{a} \dot{\bar{a}}\right) |q|^2+\left(\frac{3\bar{a}}{2}-2\bar{d}\right) |p|^2\\
+(2\bar{a}^2-2\bar{a}\bar{d}+\dot{\bar{a}}) q\T p \\
- \frac{1}{L} (1-\bar{\beta}\bar{a}) p\T (\nabla f(y)-\nabla f(q)).
\end{multline}
As in the strongly-convex case, the variable $\bar{a}(t)$ is deliberately chosen such that the cross term $q\T p$ vanishes. Direct calculations show that $\bar{a}\dot{\bar{a}} + \bar{a}^3/2$ and $3\bar{a}/2 -2\bar{d}$ vanish, which yields
\begin{multline*}
\dot{\bar{V}}\leq -\bar{a} \bar{V}
- \frac{1}{L} (1-\bar{\beta}\bar{a}) p\T (\nabla f(y)-\nabla f(q)).
\end{multline*}
The smoothness and the convexity of $f$ enables us to bound the last term, resulting in
\begin{equation}
\dot{\bar{V}}\leq \begin{cases} -\bar{a} \bar{V} - (1-\bar{\beta}\bar{a}) \bar{\beta} |p|^2, &t\in [0,1) \\
-\bar{a}\bar V, &t\geq 1.\end{cases}
\end{equation}
The case analysis is due to the fact that $\bar{\beta} < 0$ for $t\in [0,1)$ and $\bar{\beta}\geq 0$ for $t\geq 1$. 

The dynamics are uniformly Lipschitz continuous, and hence the rate of growth of the trajectories is bounded. As a result, there exists a constant $\bar{C}>0$ such that $\bar{V}(1)\leq \bar{C} (|q(0)|^2+|p(0)|^2)$. Applying Gr\"{o}nwall's inequality yields
\begin{align*}
\bar{V}(t) \leq \frac{9}{(t+2)^2} \bar{V}(1), \quad t\geq 1,
\end{align*}
and therefore,
\begin{equation*}
\frac{1}{L} f(q(t)) \leq \bar{V}(t) \leq \frac{9\bar{C}}{(t+2)^2} (|q(0)|^2+|p(0)|^2),
\end{equation*}
for all $t\geq 1$, which establishes the proposition. It can be shown that $\bar{C}\leq 5/6$. \qed

\subsection{Generalization of Prop.~\ref{Prop:Contraction} to dimensions $n\geq 1$}\label{App:Gen}
We follow the notation and definitions of \citet{RudinAnalysis} (Chapter 10), and for simplicity, we assume that the function $f$ is sufficiently smooth. We say that $E\subset \mathbb{R}^{2n}$ is a smooth two-dimensional surface with boundary if there exists a smooth map $\Phi: \mathbb{R}^2 \rightarrow \mathbb{R}^{2n}$ such that $\Phi(S)=E$, where $S$ is the unit simplex in $\mathbb{R}^2$.

\begin{proposition}
Let $\Gamma_k$ be a smooth two-dimensional surface with boundary and define
\begin{equation}
A_k:=\int_{\Gamma_k} \text{d}q_k \wedge \text{d}p_k,
\end{equation}
where $\text{d}q_{k} \wedge \text{d}p_{k}$ refers to the canonical symplectic form in $\mathbb{R}^{2n}$.
The integration scheme \eqref{eq:SympEuler1} and \eqref{eq:SympEuler2} maps $\Gamma_{k}$ to $\Gamma_{k+1}$ such that
\begin{align}
&A_{k+1}-A_{k}=-T_\text{s} \int_{\partial \Gamma_{k}} f_\text{NP}(q_k,p_k) \text{d}q_{k} \label{eq:proofeq1}\\
&=-T_\text{s} \int_{\Gamma_k} 2d + \frac{\beta}{L} \Delta f(q_k+\beta p_k) ~\text{d}q_{k}\wedge\text{d}p_{k}.\label{eq:proofeq2}
\end{align}
\end{proposition}

\begin{proof}
Let $\psi: \mathbb{R}^{2n} \rightarrow \mathbb{R}^{2n}$ denote the map induced by \eqref{eq:Euler1} and \eqref{eq:Euler2}, and let $\Phi': \mathbb{R}^2 \rightarrow \mathbb{R}^{2n}$ be a smooth map such that $\Phi'(S)=\Gamma_k$, where $S$ is the unit simplex in $\mathbb{R}^2$. Then $\Gamma_{k+1}=\psi(\Gamma_k)=\psi \circ \Phi'(S)$ is likewise a smooth two-dimensional surface with boundary.

The integration scheme \eqref{eq:Euler1} and \eqref{eq:Euler2} can be rewritten as a variational equality:
\begin{multline}
\delta q_k\T (p_{k+1}-p_{k}) -\delta p_{k+1}\T (q_{k+1}-q_k)+ T_\text{s} \delta q_k\T \nabla f(q_k) \\
+ T_\text{s}  \delta p_{k+1}\T \nabla T(p_{k+1})-T_\text{s} \delta q_k\T f_\text{NP}(q_k,p_k)=0,\label{eq:varFor}
\end{multline}
for all $\delta q_k\in \mathbb{R}^n$ and $\delta p_{k+1} \in \mathbb{R}^n$. Let the points $(q_k,p_k)$ along the boundary $\partial \Gamma_k=\Phi'(\partial S)$ be parametrized by $s\in [0,1]$, that is $(q_k(s),p_k(s))\in \partial \Gamma_k$ for all $s\in [0,1]$ and $q_k(1)=q_k(0)$, $p_k(1)=p_k(0)$. For each $s\in [0,1]$, every $(q_k(s),p_k(s))$ will be mapped to a pair $(q_{k+1}(s),p_{k+1}(s))$, which forms (by continuity of the integration scheme) the boundary of $\partial \Gamma_{k+1}$~(see, for example, \cite{RudinAnalysis} (p.~93)). We choose $\delta q_k=\text{d}q_k/\text{d}s~\text{d}s$ and $\delta p_{k+1}=\text{d}p_{k+1}/\text{d}s~ \text{d}s$ in \eqref{eq:varFor} and integrate over $s\in [0,1]$. This yields, after rearranging terms,
\begin{align}
\int_0^1& \left( -p_{k}\T \frac{\text{d} q_k}{\text{d} s} - q_{k+1}\T \frac{\text{d} p_{k+1}}{\text{d} s} - T_\text{s} f_\text{NP}\T(q_k,p_k) \frac{\text{d} q_k}{\text{d} s} \right)~\text{d}s \nonumber\\
&+\int_0^1 \text{d}\left(q_k\T p_{k+1} + T_s H(q_k,p_{k+1})\right)=0,
\end{align}
where the argument $s$ has been omitted to shorten notation.
The second part is a total differential and vanishes due to the fact that $\partial \Gamma_k$ and $\partial \Gamma_{k+1}$ are closed curves. This implies that 
\begin{equation}
\int_{\partial \Gamma_k} \hspace{-4pt}p_k \text{d} q_k + \int_{\partial \Gamma_{k+1}} \hspace{-12pt}q_{k+1} \text{d} p_{k+1} = -T_\text{s} \int_{\partial \Gamma_{k}} f_\text{NP}(q_k,p_k) \text{d}q_k.
\end{equation}
Applying Stokes' theorem \cite{RudinAnalysis}(p.~272) to the left-hand side yields \eqref{eq:proofeq1}. The formula \eqref{eq:proofeq2} results from applying Stokes' theorem to the right-hand side. \qed
\end{proof}

\subsection{Proof of Prop.~\ref{Prop:ContractionIII}}\label{App:ProofContractionIII}
Let $\partial \Gamma(0)$ be the positively oriented contour along the boundary of an energy level set, and consider its evolution through the continuous-time dynamics, $\partial \Gamma(t)$, $t\in [0,T_\text{s}]$, and the discretized dynamics, $\partial \Gamma_1$. By the mean value theorem it follows from \eqref{eq:ctac} that there exists a time instant $t_\text{m}\in (0,T_\text{s})$ such that 
\begin{equation}
A(T_\text{s})=A(0) - T_\text{s} \int_{\Gamma(t_\text{m})}2d +\frac{\beta}{L} \Delta f(q+\beta p)~\text{d}q \text{d}p,
\end{equation}
where $A(t)$ refers to the area enclosed by $\partial \Gamma(t)$. Due to the fact that energy is dissipated as the continuous-time system evolves, it holds that $\Gamma(t)\subset \Gamma(0)$ for all $t>0$. Therefore it follows that 
\begin{align}
A(T_\text{s})&\geq A(0)-T_\text{s} \int_{\Gamma(0)}2d +\frac{\beta}{L} \Delta f(q+\beta p)~\text{d}q \text{d}p\\
&=A_{1},
\end{align}
where $A_{1}$ denotes the (signed) phase-space area enclosed by $\partial \Gamma_{1}$. \qed

\subsection{Proof of Prop.~\ref{Prop:Hom}}\label{App:ProofOfHom}
The map $\psi$, induced by \eqref{eq:Euler1} and \eqref{eq:Euler2} is Lipschitz continuous due to the fact that $\nabla f$ is Lipschitz continuous. As mentioned earlier, $\psi$ can be divided into two parts, a contraction step and a semi-implicit Euler step based on the conservative part of the system. We show that both steps are bijective, the result then follows by the invariance of domain theorem \cite{Munkres}. 

We start with the contraction step 
\begin{equation}
\bar{p}_{k+1}=p_k+T_\text{s} (-2d p_k -\frac{1}{L}(\nabla f(q_k+\beta p_k)-\nabla f(q_k))), \label{eq:contrStep1}
\end{equation}
and apply fixed-point iteration to solve for $p_k$, given $q_k=\bar{q}_{k+1}$ and $\bar{p}_{k+1}$. Let $p_k^j$ be the $j$th iterate, which is updated according to
\begin{align*}
p_k^{j+1}&=\bar{p}_{k+1}+T_\text{s} (2d p_k^j +\frac{1}{L}(\nabla f(q_k+\beta p_k^j)-\nabla f(q_k)))\\
&=:F_{\bar{p}_{k+1}, q_k}(p_k^j),
\end{align*}
with $p_k^0=0$.
Note that $F_{\bar{p}_{k+1}, q_k}$ is Lipschitz continuous, as for any $p_1, p_2\in \mathbb{R}^n$,
\begin{align*}
|F_{\bar{p}_{k+1}, q_k}(p_1)-F_{\bar{p}_{k+1}, q_k}(p_2)|&=|T_\text{s} (2d(p_1-p_2)+\\
&\hspace{-2cm}\frac{1}{L}(\nabla f(q_k +\beta p_1)-\nabla f(q_k + \beta p_2)))|\\
&\leq T_\text{s} (2d+\beta) |p_1-p_2|.
\end{align*}
Due to the fact that $2d+\beta=1$, the map $F$ is a contraction for any $T_\text{s} \in (0,1)$ and the fixed-point iteration is guaranteed to converge. Thus, $p_k=\lim_{j\rightarrow \infty} p_k^j$ exists and is unique, which shows that the contraction step is a bijection.

The semi-implicit Euler step according to  \eqref{eq:SympEuler1} and \eqref{eq:SympEuler2} can be inverted in the following way:
\begin{align}
\bar{q}_{k+1}&=q_{k+1}- T_\text{s} \nabla T (p_{k+1})\\
\bar{p}_{k+1}&=p_{k+1}+ T_\text{s} \frac{1}{L} \nabla f(\bar{q}_{k+1}) \nonumber\\
&=p_{k+1}+ T_\text{s} \frac{1}{L} \nabla f(q_{k+1}- T_\text{s} \nabla T (p_{k+1})),
\end{align}
which maps $(q_{k+1},p_{k+1})$ to $(\bar{q}_{k+1}, \bar{p}_{k+1})$. \qed

\subsection{Proof of Prop.~\ref{Prop:HardConvergence}}\label{App:HardConvergence}
The proof is inspired by \citet{NesterovBook}(Section 2.2.1). 
Without loss of generality we set the Lipschitz constant $L$ to one and introduce the following change of coordinates $(q_k,p_k) \rightarrow (\hat{q}_k,\hat{p}_k)$:
\begin{equation}
    \hat{q}_{k}:=q_k + \left(\frac{\beta}{1-2dT_\text{s}}-T_\text{s}\right) p_k, \quad
    \hat{p}_{k}:=p_k.
\end{equation}
The discrete-time dynamics expressed in the $(\hat{q}_k,\hat{p}_k)$-coordinates read as
\begin{align}
    \hat{q}_{k+1}&=y_k - T_\text{s} \tau \nabla f(y_k),\\
    \hat{p}_{k+1}&=(1-2dT_\text{s}) \hat{p_k} -T_\text{s} \nabla f(y_k),
\end{align}
where $\tau:=\beta/(1-2dT_\text{s})$ and 
\begin{align}
    y_k&=q_k + \beta p_k\\
    &=\hat{q}_k + \underbrace{\left(T_\text{s}-\frac{2d\beta T_\text{s}}{1-2dT_\text{s}}\right)}_{=:\hat{\beta}} \hat{p}_k.
\end{align}
The coordinate transformation is motivated by the fact that the smoothness and convexity of the objective function can be used to conclude $f(\hat{q}_{k+1}) \leq f(y_k)$ (see also \citet{NesterovBook}(p.~76, General scheme of optimal method, step 1.c).)),
\begin{multline}
f(\hat{q}_{k+1})-f(y_k) \leq 
- |\nabla f(y_k)|^2 \underbrace{\left(T_\text{s} \tau - \frac{T_\text{s}^2 \tau^2}{2}\right)}_{\geq 0, \forall T_\text{s}\in (0,1], \kappa\geq 1}.\label{eq:proofeqtmp1}
\end{multline}
We claim that the following function is a Lyapunov function
\begin{multline}
    \hat{V}(\hat{q}_k,\hat{p}_k)=\frac{1}{2} \Big|d \hat{q}_k + (1-d\tau) \hat{p}_k\Big|^2 + f(\hat{q}_k).
\end{multline}
It can be checked that the function is indeed positive definite for all $T_\text{s} \in (0,1]$ and all $\kappa\geq 1$. Moreover, the evolution of $\hat{V}$ along the discrete-time trajectories is given by
\begin{align*}
    \hat{V}_{k+1}-\hat{V}_{k} &\leq -d T_\text{s} \hat{q}_k\T \nabla f(y_k) -\hat{\beta} \hat{p}_k\T \nabla f(y_k)\\
    &+C_{\text{pn}} \hat{p}_k\T \nabla f(y_k) + C_{\text{pp}} |\hat{p}_k|^2 + C_{\text{pq}} \hat{p}_k\T \hat{q}_k\\
    &+C_{\text{nn}} |\nabla f(y_k)|^2 + f(y_k)-f(\hat{q}_k),
\end{align*}
where \eqref{eq:proofeqtmp1} has been used to relate $f(\hat{q}_{k+1})$ to $f(y_k)$, $C_{\text{pp}},C_\text{pn},C_\text{pq},C_\text{nn}$ are constants (dependent on $\kappa$ and $T_\text{s}$), $\hat{V}_{k}$ denotes $\hat{V}(\hat{q}_k,\hat{p}_{k})$, and $\hat{V}_{k+1}$ denotes $\hat{V}(\hat{q}_{k+1},\hat{p}_{k+1})$.
It follows from the strong convexity and the smoothness of $f$ that
\begin{align*}
-y_k \nabla f(y_k) + \frac{1}{2\kappa} |y_k|^2 + f(y_k) &\leq 0,\\
f(y_k)-\hat{\beta} \hat{p}_k\T \nabla f(y_k) + \frac{1}{2\kappa} \hat{\beta}^2 |\hat{p}_k|^2 &\leq f(\hat{q}_k),
\end{align*}
which can be used to conclude
\begin{align*}
&-d T_\text{s} \hat{q}_k\T \nabla f(y_k) - \hat{\beta} \hat{p}_k\T \nabla f(y_k) + f(y_k)-f(\hat{q}_k)\\
&\leq -d T_\text{s} f(\hat{q}_k) - \frac{d T_\text{s}}{2 \kappa} |\hat{q}_k|^2 - \frac{d \hat{\beta} T_\text{s}}{\kappa} \hat{q}_k\T \hat{p}_k - \frac{1}{2\kappa} \hat{\beta}^2 |\hat{p}_k|^2.
\end{align*}
As a result, the difference $\hat{V}_{k+1}-\hat{V}_{k}$ is upper bounded by
\begin{align*}
    \hat{V}_{k+1}-\hat{V}_{k} &\leq - d T_\text{s} \hat{V}_k + C_{\text{pn}} \nabla f(y_k)\T \hat{p}_k\\
    &+\left( C_\text{pp} - \frac{\hat{\beta}^2}{2\kappa} + \frac{d T_\text{s}}{2} (1-d\tau)^2 \right) |\hat{p}_k|^2 \\
    &+\left( C_\text{pq} - \frac{d T_\text{s} \hat{\beta}}{\kappa} + d^2 T_\text{s} (1-d\tau) \right) \hat{p}_k\T \hat{q}_k\\
    &+d T_\text{s} \left(\frac{d^2}{2}-\frac{1}{2\kappa}\right) |\hat{q}_{k}|^2 + C_\text{nn} |\nabla f(y_k)|^2.
\end{align*}
Rearranging terms we obtain
\begin{multline}
\hat{V}_{k+1}-\hat{V}_{k} \leq -d T_\text{s} \hat{V}_k + \\
(\hat{p}_k\T, \nabla f(y_k)\T, \hat{q}_k\T) (M \otimes I) (\hat{p}_k\T, \nabla f(y_k)\T, \hat{q}_k\T)\T,
\end{multline}
where the matrix $M\in \mathbb{R}^{n\times n}$ is constant (depends on $\kappa$ and $T_\text{s}$), $\otimes$ denotes the Kronecker product, and $I\in \mathbb{R}^{n\times n}$ the identity. It can be shown, for instance by applying Sylvester's criterion combined with Taylor expansions, that for every $T_\text{s}\in (0,1)$, the matrix $M$ is negative definite for $\kappa \rightarrow \infty$. For $T_\text{s}=1$, the matrix $M$ is negative semi-definite for all $\kappa\geq 1$.\footnote{In fact, numerical evaluations indicate that $M$ is negative semi-definite for all $\kappa\geq 3$ and all $T_\text{s} \in (0,1]$.} Thus, for every $T_\text{s} \in (0,1]$, we obtain, for sufficiently large $\kappa$,
\begin{equation}
\hat{V}_{k+1}-\hat{V}_{k} \leq -d T_\text{s} \hat{V}_k,
\end{equation}
yielding the desired result. \qed

\end{document}